\documentstyle[12pt]{article}\textwidth 170mm\textheight 220mm\topmargin -2cm
\begin{document}
\bibliographystyle{unsrt}\begin{center}{\large\bf New Coordinate Systems and
Potentials on a Plane Providing Unusual Intergral of Motion}\\[5mm]
Z.~Ya.~Turakulov\\{\it Inter-University Centre for Astronomy and Astrophysics,
Post Bag 4,\\ Ganeshkhind, Pune~411~007, India\\ and Institute of Nuclear
Physics\\ Ulugbek, Tashkent 702132, Uzbekistan\\ (e-mail:zunn@tps.uz)}
\end{center}\begin{abstract}

We discuss a new class of coordinate systems for a plane, which provide an 
analytical representation of arbitrary straightline, and then define the form
of potential on the plane, under which the equations of motion of a mass point are
solvable. In both cases we use the Hamiltonian approach, reducing the problem to the
Hamilton-Jacobi equation. The coordinate systems in question are defined in a 
somewhat unusual way, and, correspondingly, the potentials we obtain, also
possess some unusual properties. Since the equations of motion are solvable within
this potential, the second integral of motion is defined too, and, like
the coordinate system and potential, it possesses some nonstandard properties.
\end{abstract}\vspace{1cm}\section{Introduction}

All known integrals of motion of conservative mechanical systems can be
obtained as constants appearing as the result of variables separation in the
Hamilton-Jacobi equation for given system. They always are either linear or
quadric functions of generalized velocities of the system; the earlier are
provided by its symmetry while the latter appear when system has no symmetry
group. No other forms of integrals of motion are known. Therefore it is
interesting to explore possibility of existence of systems possessing integrals
of motion of other form.

In this paper we make an attempt to find out systems with other integrals of
motion in particular case of two degrees of freedom. The systems under
consideration describe motion of particle on plane under action of given
potential. The work is organized as follows. As new class of desired systems is
found we pass to the next stage and search for potentials which can be included
into the equation without losing separability of the equation. Finally we
discuss properties of the potentials satisfying these conditions.

To show what do we mean by new class of coordinate systems, let us start with
the commonplace definition of coordinate transformation. This definition reads
that coordinate transformation $\{x^i\}\rightarrow\{y^a\}$ is just introducing
some functions $y^a$ of coordinates\begin{equation} y^a=y^a(x^i)
\end{equation} which satisfy well-known conditions, and it suffices to know
that the inverse transformation\begin{equation} x^i=x^i(y^a)\end{equation}
exists. From practical point of view this is not sufficient, because, besides
definition of new coordinates (1) one needs explicit form of metric in the new
system. Transformation of the metric requires that the equations (1) are solved
analytically such that explicit form of old coordinates $x^i$ as functions of
new ones is known. In other words, the difference between commonplace
definition and practical needs is that the earlier assumes only existence of
the functions $x^i(y^a)$ and the latter require them to be given explicitly.
This means that, in fact, very few kinds of coordinate systems can be used
practically, because analytical solvability of the equations (1) means that
they actually can be reduced to algebraic equations of not more than second
order, otherwise they are not solvable. However, all coordinate systems of this
sort are known and solvability of Hamilton-Jacobi equation has been studied in
all of them. Therefore, in order to find out new usable coordinate systems we
try to broaden the scope and consider systems which cannot be defined
analytically as in the equation (1). By new class of corrdinate systems we mean
systems defined such an unusual way.

There is strong interdependence between solvability of the equations (1) and
that of free motion of a particle. If, for example, $\{x^1\}$ are Cartesian
coordinates, solutions of the equations (1) specify strightlines which, at the
same time, are trajectories of free motion of a particle. Therefore, if
equations (1) cannot be solved analytically, solvability of equations of
dynamics is doubtful because, otherwise, Cartesian coordinates can be restored
from straightlines obtained as trajectories. On the other hand, solvability
means possibility to express solution in terms of known objects which does not
depend on approach to the problem. Therefore in this work solvability of both
equations (1) and equations of dynamics are considered to be the same.
\section{Strange coordinate systems and their metrics}

All coordinate transformations used till now specify new coordinates $y^a$
analytically as certain functions of the old ones. Another possibility is to
define new coordinates geometrically. To do it specfy two foliations of plane
or its domain with two families of curves and introduce two numerical
parameters each of which labels curves of one family. These two parameters are
new coordinates which, however, are not defined as functions of old ones,
moreover, when defining them we did not need any `old' coordinates to be
introduced previously. In fact, we use here the natural way of introducing
coordinates on plane for the first time, as if no coordinates on plane have
been introduced before. As it is done, our task is to obtain metric of this
system from geometric properties of the foliations and definition of the
parameters introduced.

Let us now pass to practical implementation of this program. Consider a simple
arch of curve on plane. By simplisity we mean that it is smooth and has no
intersections with straightlines tangent to it. Hereafter we call it the basic
curve of the coordinate system. Consider now all straight rays tangent to the
curve and starting from each its point in (locally) one direction. These
straight rays constitute foliation of a domain on plane, whose boundary
consists of the basic curve and two rays starting from its endpoints. Now we
build a coordinate system for this domain using the rays as coordinate lines.
For this end we label each ray with the value of angle $\psi$ between the ray
and some fixed direction on the plane. Thus, $\psi$ is one of new coordinates.
To introduce the second coordinate, note that all rays $\psi=const$ are
orthogonal to all evolvents of the curve all rays are tangent to, hence, in
order to obtain an orthogonal coordinate system we must use foliation of the
domain with the evolvents. Now it remains to parametrize the family of
evolvents with a numerical parameter and find out the metric of coordinate
system defined this way.

Evolvents of a curve on plane \cite{1,2} are known to be equidistant lines
(parallels) \cite{3}. Therefore it is natural to label them with the length
parameter $R$. If this parameter is defined properly it specifies a function
on the domain, which satisfies the Hamilton-Jacobi equation for straightlines.
Gradients of $R$ and $\psi$ as functions on plane are orthogonal and radius of
curvature of evolvent $R=const$ equals $R-l(\psi)$ where $l(\psi)$ depends
only on shape of the basic curve. Square gradient of the function $\psi$ is
inverse square of radiusof curvature:\begin{equation} <dR,dR>=1,\ <dR,\psi>=0,\
<d\psi,d\psi>=[R-l(\psi)]^{-1}.\end{equation} These equalities specify the
metric of the coordinate system$\{R,\psi\}$. \section{Arbitrary straightlines
and Cartesian coordinates}

Each solution of Hamilton-Jacobi equation specifies one congruence of
straightlines. We use this fact to obtain analytical representation of
arbitrary straightline in the coordinates $\{R.\psi\}$. In this case the
equation has the form\begin{equation}<dS,dS>=1,\end{equation} and, due to the
equations (3) can be written as$$\left(\frac{\partial S}{\partial R}\right)^2+
\frac{1}{[R-l(\psi)]^2}\left(\frac{\partial S}{\partial\psi}\right)^2=1.$$In
order to to separate variables in Hamilton-Jacobi equation one assumes that the
function to be found has the form of sum of single variable functions each of
which depends only on one coordinate. In this case this means that the desired 
function must be taken in the form $S(R,\psi)=f(R)+g(\psi)$. Note that this is
the only known method of separating variables in this equation. As the equation
just obtained cannot be separated this way, to solve it we have to think out a 
new ansatz about the desired function, because the standard one does not work
in this case. Nevertheless, it will be shown below that this equation reduces
to some ordinary differential equation, and obtain the form of the desired
function under which variables separate in metrics like (6). This form will be 
used when separating variables in more general case of particle motion in some
potential.

We reduce the equation (4) to an ordinary differential equation in two steps.
First, we introduce a 1-form $\pi$ of unit norm, and second, we require that
this form is closed:\begin{equation}<\pi,\pi>=1,\ d\pi=0.\end{equation} Then,
as such a 1-form is obtained, we put $dS=\pi$ that gives solution of the
equation (4). To obtain a norm one 1-form we introduce an orthonormal frame of
1-forms $\{\nu^a\}$. Its form follows from the equations (3):\begin{equation}
\nu^1=dR,\ \nu^2=[R-l(\psi)]d\psi,\end{equation}and exterior derivatives of
$\nu^a$'s are\begin{equation} d\nu^1=0,\ d\nu^2=dR\wedge d\psi
\equiv[R-l(\psi)]^{-1}\nu^1\wedge\nu^2.\end{equation}

Now, arbitrary norm one 1-form can be represented as follows:\begin{equation}
\pi=\nu^1\sin f+\nu^2\cos f\end{equation} where $f$ is an arbitrary function on
the plane. It remains to find it out from the secone equation (5):$$0=d\pi=
-\cos f\nu^1\wedge df+\sin fdf\wedge\nu^2+\cos fd\nu^2.$$ Substituting  the
equations (6) and (7) gives:
$$\{-[R-l(\psi)]^{-1}\cos ff_\psi-\sin ff_R+[R-l(\psi)]^{-1}\cos f\}\nu^1\wedge
\nu^2=0.$$ One particular solution of this equation is evident: $f=\psi+const$.
However, since it contains an arbitrary constant no more general solution is
needed,and we put$$dS=\sin(\psi-\psi_0)dR+\cos(\psi-\psi_0)[R-l(\psi)]d\psi.$$
The desired solution is$$S=R\sin(\psi-\psi_0)-\int l(\psi)\cos(\psi-\psi_0).$$
Note that differentiating this function on the constant $\psi_0$ is equivalesnt
to change of the constant. Due to the Jacobi theorem \cite{4} this means that
lines $S=const$ are straightlines, consequently the function S specifies one of
Cartesian coordinates or their linear combination. Alternatively, Cartesian
coordinates can be introduced as two solutions of the Hamilton-Jacobi equation
$x^1$ and $x^2$ with orthogonal gradients $<dx^1,dx^2>=0$:\begin{equation}
x^1=R\sin\psi-\int l(\psi)\cos\psi d\psi,\ x^2=R\sin\psi+\int l(\psi)\sin\psi
d\psi.\end{equation} It is seen now that the coordinares $R$ and $\psi$ cannot
be expressed analytically as functions of Cartesian coordinates (9), hence,
equation (2) has purely formal meaning. \section{Inclusion of potential}

As seen from the result obtained, we have found how to separate Hamilton-Jacobi
equation in metrics of the form (3). For this end the function to be found is
to be taken in the form\begin{equation} S=Rf(\psi)+g(\psi).\end{equation} As
this is clear now, let us look for possible form of force function which can be
included into the Hamilton-Jacobi equation without losing separability.
Hamilton-Jacobi equation for particle of mass $m$ and energy $E$, moving in
potential $V(R,\psi)$ has the form\begin{equation}
\frac{1}{2m}<dS,dS>+V(R,\psi)=E.\end{equation} Assuming that the function to be
found has the form (10), we have$$dS=f(\psi)\nu^1+
\frac{Rf'+g'}{R-l(\psi)}\nu^2$$ and it is seen that condition of separability
reads:\begin{equation} g(\psi)=-\int l(\psi)f'(\psi)d\psi\end{equation} such
that $dS=f(\psi)\nu^1+f'(\psi)\nu^2.$ Then the potential is function of single
variable $\psi$, and if we denote $2m[E-V(\psi)\equiv p^2(\psi)$ the
Hamilton-Jacobi equation (11) reduces to the following ordinary differential
equation:\begin{equation} f'^2+f^2=p^2.\end{equation}Finally, solution of the
equation (11) appears as$$S=Rf(\psi,C)-\int l(\psi)f'(\psi,C)d\psi$$where $C$
is arbitrary constant of general solution of the equation (13). Note that this
constant plays the role of the second integral of motion. All known integrals
of motion can be obtained from the procedure of variables separation in
Hamilton-Jacobi equation and are linear or quadratic functions of the particle
velocity. This one appears later, after variables are separated and
Hamilton-Jacobi equation s reduced to an ordinary differential equation. It is
seen that this integral of motion is not linear or quadratic on velocity,
moreover, probablly, it cannot be expressed in terms of known functions of
coordinates and velocity components.\section{Conclusion}

Thus, we have constructed a class of strange coordinate systems $\{R,\psi\}$
which cannot be expressed in terms of known functions of any standard
coordinates on plane. However, Cartesian coordinates can be expressed
analytically as functions of $R$ and $\psi$ and this expression is given by the
equations (9). The metrics of all coordinate systems of this class has one and
the same form (3) and differ only in the form of the function $l(\psi)$
specified by the basic curve of the system. Hamilton-Jacobi equation for
straightlines on the plane separates in these metrics, and there exists a class
of potentials in which Hamilton-Jacobi equation for a mass point also
separates. All potentials of this class depend on only one coordinate $\psi$
and also possess some strange properties. Like the coordinates, these
potentials cannot be expressed in terms of known functions of any standard
coordinates on plane, and, since equipotential lines are rays $\psi=const$,
which have envelope, the potentials are non-univalent on the basic curve. The
second integral of motion of particle moving in potential of this sort is not
of known form, at least, it is certainly not of first or second order on
velocity components, and seems to belong to the same class of unknown
functions of coordinates and velocities, as the potentials.
{\bf Acknowledgment:}The author thanks the Third World Academy of Sciences for
financial support and IUCAA for warm hospitality which made this work possible.
\end{document}